\documentclass[12pt,eps]{article}
\usepackage{hadronic}

\usepackage{graphicx}
\usepackage{epstopdf}
\newcommand{\ed}{\end{document}}
\newcommand{\be}{\begin{equation}}
\newcommand{\ee}{\end{equation}}
\newcommand{\bc}{\begin{center}}
\newcommand{\ec}{\end{center}}
\newcommand{\ba}{\begin{array}}
\newcommand{\ea}{\end{array}}

\newcommand{\cl}{\centerline}

\begin{document}
\thispagestyle{empty}

{\cl {\Huge \bf 
GENERALIZED ALEXANDER  }
\vspace{0.4cm}
{\cl {\Huge \bf 
POLYNOMIAL INVARIANTS\footnote{Based on invited talks given at the 5th Petrov International Symposium on {\it High Energy Physics, Cosmology
and Gravity}, Kyiv (Ukraine), April 29--June 15, 2012 and partially supported by 
 the Project No. 1202.094-12 of the
Central European Initiative Cooperation Fund.
}}
\vspace{1cm}
{\cl{
{Anatoliy M. PAVLYUK}}}
\vspace*{0.2cm}
\begin{center} {
Bogolyubov Institute for
Theoretical Physics}\\
\vspace*{0.2cm}
{
National Academy of Sciences of Ukraine}\\
\vspace*{0.2cm}
{
 Metrologichna str. 14b, Kyiv 03143, Ukraine}\\
 \vspace*{0.2cm}
{
e-mail: pavlyuk@bitp.kiev.ua}\\
\end{center}

\vspace{1cm}
\cl\Large{ {\bf  Abstract}}
We propose an algorithm which allows to derive 
the generalized Alexander polynomial invariants
of knots and links
with the help of the $q,p-$numbers, appearing in 
bosonic two-parameter quantum algebra.
These polynomials turn into HOMFLY ones by applying
special parametrization. The Jones polynomials can be also obtained
by using this algorithm.

\vspace{0.3cm}
\noindent
{\bf Keywords:} knots and links; skein relationship; recurrence relation; $q-$numbers;
$q,p-$numbers; polynomial invariants of knots and links; Alexander polynomials;
generalized Alexander polynomials;
HOMFLY polynomials; Jones polynomials.

\vspace{0.3cm}
\noindent
{\bf Mathematics Subject Classification 2010:}  57M27,  17B37,  81R50.

\newpage

\noindent
\section{
Introduction}

The aim of this paper is to generalize one-parameter Alexander polynomial invariants,
one of the main characteristics of knots and links,
to two-parameter generalized Alexander polynomial invariants.

 First, we recall some basic notions
of the knot theory. 
Applying to an initial link (knot) $L_{+}$ so called ``surgery'' \  operation  - elimination of a 
crossing - we obtain a simpler link/knot
$L_{O}$. Applying to the same initial link (knot) $L_{+}$ another ``surgery'' \ operation - 
 switching of a crossing -  
we obtain another simpler link/knot
$L_{-}$. 

It is postulated:\\
1) every knot and link is described by the definite polynomial;\\
2) three concrete polynomials, namely 
$P_{L_{+}}(t),\, P_{L_{O}}(t),$  $P_{L_{-}}(t)$ 
are connected with the help of
the following (geometro-algebraic) recurrence relation, which is called
the skein relationship:
 \be\label{skein}
P_{L_{+}}(t)=l_{1}P_{L_{O}}(t)+l_{2}P_{L_{-}}(t)\ee
where $l_{1}, l_{2}$ are coefficients;\\
3) the normalization condition for the unknot: 
\be\label{unknot} P_{unknot}=1\,.\ee

Applying the operation of elimination for torus knots and links $L_{n, 2}$ turns it into $L_{n-1, 2}$, 
and the switching operation turns it into $L_{n-2, 2}$, where $n$ is a positive integer number.
From these considerations and from~(\ref{skein}) it follows the following
recurrence relation :
\[ P_{L_{n+1, 2}}(t)=l_{1}P_{L_{n, 2}}(t)+l_{2}P_{L_{n-1, 2}}(t)\,,\]
or in the simpler notations:
 \be\label{skein1}
P_{n+1, 2}(t)=l_{1}P_{n, 2}(t)+l_{2}P_{n-1, 2}(t)\,.
 \ee
Thus, the form of the recurrence relation~(\ref{skein1}) for torus knots and links $L_{n, 2}$ coincides
with the form of the skein relationship~(\ref{skein}).

Recurrence relation only for torus knots $T(2m+1, 2)$ (or only for torus links $L(2m, 2)$),
where $m=0,1,2,\ldots\,.$ 
looks as follows:
 \be\label{skein2}
P_{n+2, 2}(t)=k_{1}P_{n, 2}(t)+k_{2}P_{n-2, 2}(t)\,,
 \ee
where
\be\label{k}
k_{1}=l_{1}^{2}+2l_{2},\quad k_{2}=-l_{2}^{2}.\ee
We also have
\be\label{32}
P_{1,2}=1,\quad P_{3,2}=k_{1}+k_{2}.
\ee

\vspace{0,5cm}
\noindent

\section{
Alexander polynomials}

 The Alexander polynomials $\Delta(t)$ of knots and links~\cite{Al} can be defined
by the skein relationship
 \be\label{alex-skein}
\Delta_{+}(t)-\Delta_{-}(t)=(t^{1\over2}-t^{-{1\over2}})\Delta_{O}(t),\quad \Delta_{unknot}=1.
 \ee
From (\ref{alex-skein}) (in analogy to~(\ref{skein1})) it follows the following recurrence relation
for torus knots and links $L_{n,2}(t)$:
 \be\label{alex-skein1}
\Delta_{n+1, 2}(t)=(t^{1\over2}-t^{-{1\over2}})\Delta_{n, 2}(t)+\Delta_{n-1, 2}(t)\,.
 \ee
From (\ref{alex-skein1}) (in analogy to~(\ref{skein2})) one obtains the  recurrence relation
only for torus knots $T(2m+1, 2)$ (or for torus links $L(2m, 2)$)
 \be\label{alex-skein2}
\Delta_{n+2, 2}(t)=(t+t^{-1})\Delta_{n, 2}(t)-\Delta_{n-2, 2}(t)\,.
 \ee

The Alexander polynomials of torus knots $T(n, 2)$ can be expressed through
$q$-numbers characteristic to Biedenharn-Macfarlane quantum bosonic oscillator.
The bosonic $q$-number corresponding to an integer $n$
 is defined as~\cite{Bi,Ma}
 \be \label{q-def} [n]_{q}={{\textstyle
{q^{n}-q^{-n}}}\over{\textstyle {q-q^{-1}}}}\,,
 \ee
  where $q$ is a parameter.
Some of the $q$-numbers are:
\[ 
[1]_{q}=1\,,\quad  [2]_{q}=q+q^{-1}\,,\quad
 [3]_{q}=q^{2}+1+q^{-2}\,, \quad  [4]_{q}=q^{3}+q+q^{-1}+q^{-3}, \ldots  \,.
 \]
 The recurrence relation for~(\ref{q-def})
looks as 
\be\label{q-rec}
[n+1]_{q}=(q+q^{-1})[n]_{q}-[n-1]_{q}\,.\ee

It was found that~\cite{Ga,GP2}:
 \be\label{A2m} \Delta_{2m+1,2}(t)=[m+1]_{t}-[m]_{t}\,,
 \hspace{10mm} t\equiv q \ ,
 \ee
or, since $n=2m+1$, as
 \be\label{A2n}
\Delta_{n,2}(t)=\Bigl[{{n+1}\over2}\Bigr]_{t}-\Bigl[{{n-1}\over2}\Bigr]_{t}\,.\ee

In the following section we generalize these results with the help
of $q-$numbers.

\vspace{0,5cm}
\noindent
\section{
Algorithm of obtaining of Alexander \\ polynomials from bosonic $q-$numbers}

Analyzing the results of previous sections we can formulate
an algorithm of obtaining of the Alexander skein relationship~(\ref{alex-skein}).
Afterwards this procedure will be used for obtaining another skein relations.

First step:  we introduce  polynomials 
$A_{n, 2}(q),$ which refer to torus knots  $T(2m+1,2)$,  satisfying following recurrence relation (repeating~(\ref{q-rec})):
\be\label{Askein2}
 A_{n+2, 2}(q)=(q+q^{-1})A_{n, 2}(q)-A_{n-2, 2}(q)\,.\ee
 According to~(\ref{32}):
\be\label{A32}
A_{1,2}(q)=1,\quad A_{3,2}(q)=q-1+q^{-1}.
\ee

Second step: we formulate full recurrence relation for all polynomials
$A_{n, 2}(q)$ and, thus, find corresponding skein
relationship.  
From (\ref{Askein2}) we have
$k_{1}=q+q^{-1},\,  k_{2}=-1.$
Because of (\ref{k}), 
we find 
\be\label{Al}
 l_{1}=q^{1\over2}+q^{-{1\over2}}, \quad l_{2}=1\,.
\ee
Therefore
\be\label{Aalex-skein1}
A_{n+1, 2}(q)=(q^{1\over2}-q^{-{1\over2}})A_{n, 2}(q)+A_{n-1, 2}(q)\,.
 \ee
From~(\ref{Aalex-skein1}) (in anology with~(\ref{skein}) and (\ref{skein1}))
we obtain the following skein relationship:
\be\label{Aalex-skein}
A_{+}(q)-A_{-}(q)=(q^{1\over2}-q^{-{1\over2}})A_{O}(q).
 \ee

Third step: we find an expression for torus knots $A_{2m+1,2}(q)$.
In analogy with~(\ref{AA2m}), we put
\be\label{AA2m} A_{2m+1,2}(q)=a_{1}(q)[m+1]_{q}-a_{2}(q)[m]_{t}\,,
 \hspace{10mm} t\equiv q \ ,
 \ee
Using (\ref{q-def}), (\ref{A32}) and (\ref{AA2m}), we find 
$a_{1}(q)=1,\ a_{2}(q)=1.$
Therefore,
\be\label{AA2m-} A_{2m+1,2}(q)=[m+1]_{q}-[m]_{q}\,.
 \ee

In general, we described three-step procedure of obtaining of: 1) skein relationship
of knots and links, and 2) expression for polynomial invariants of torus knots
$T(2m+1,2),$ from structural functions of bosonic deformed oscillators. 
In particular, we obtained the formulas~(\ref{Aalex-skein}),~(\ref{AA2m-}),
which coincides with those 
for the Alexander polynomial invariants~(\ref{alex-skein}),~(\ref{A2m}) (if $q\equiv t$).

\vspace{0,5cm}
\noindent

\section{
Generalized Alexander polynomials $A(q,p)$ from $q,p$-numbers}

In this section we use the proposed three-step algorithm to obtain the
generalized Alexander polynomials $A(q,p)$ from $q,p$-numbers,
which reduce to the Alexander polynomials if $p=q^{-1}.$

The $q,p$-number corresponding to integer number $n$ is introduced
as~\cite{CJ} 
 \be \label{qp-def}
[n]_{q,p}={{\textstyle {q^{n}-p^{n}}}\over{\textstyle {q-p}}}\,, \ee
where $q\,,p$ are some complex parameters. If $p=q^{-1}\,,$ then
$[n]_{q,p}=[n]_{q}\,.$
Here are some of the $q,p$-numbers:
$$ [1]_{q,p}{=}1\,,\ [2]_{q,p}{=}q+p\,,\  [3]_{q,p}{=}q^{2}+qp+p^{2}\,,\  
 [4]_{q,p}{=}q^{3}+q^{2}p+qp^{2}+p^{3}, \ldots\ .$$
The recurrence relation for $q,p$-numbers
is  
\be\label{qp-rec}
[n+1]_{q,p}=(q+p)[n]_{q,p}-qp[n-1]_{q,p}\,.\ee

First, in anology with previous section,
on the base of (\ref{qp-rec}) we introduce
 polynomials $A_{n,2}(q,p)$, which generalize the Alexander polynomials:  
\be\label{A-knot}
 A_{n+2, 2}(q,p)=(q+p)A_{n, 2}(q,p)-qpA_{n-2, 2}(q,p)\,.\ee
Thus from normalization condition and (\ref{32})
\be\label{gA-32}
A_{1,2}(q,p)=1,\quad A_{3,2}(q,p)=q-qp+p.\ee

Second, from (\ref{A-knot}) it also follows
\[ k_{1}=l_{1}^{2}+2l_{2}=q+p, \quad k_{2}=-l_{2}^{2}=-qp\,.\]
From here one finds
\[ l_{2}=q^{1\over2}p^{1\over2}, \quad l_{1}=q^{1\over2}-p^{1\over2}\,,\]
which leads to the generalized Alexander skein relationship~\cite{GP1}:
\be\label{gA}
 A_{+}(q,p)=(q^{1\over2}-p^{1\over2})A_{O}(q,p)+q^{1\over2}p^{1\over2}A_{-}(q,p).\ee
Formula (\ref{gA}) can be written in the form
\be\label{gA-}
q^{-{1\over4}}p^{-{1\over4}} A_{+}(q,p)-q^{1\over4}p^{1\over4}A_{-}(q,p)=
(q^{1\over4}p^{-{1\over4}}-q^{-{1\over4}}p^{1\over4})A_{O}(q,p)\ee
By putting $p=q^{-1},$ the generalized Alexander skein relationship  
turns into the Alexander skein relationship~(\ref{alex-skein}).

Third,
we take
\be\label{AA2m-} A_{2m+1,2}(q,p)=a_{1}(q,p)[m+1]_{q,p}-a_{2}(q,p)[m]_{q,p}\,.
 \ee
From  (\ref{gA-32})  we have 
$a_{1}(q,p)=1,\ a_{2}(q,p)=qp.$
Therefore,
\be\label{AA2m-} A_{2m+1,2}(q,p)=[m+1]_{q,p}-qp[m]_{q,p}\,.
 \ee

\vspace{0,5cm}
\noindent

\section{
Generalized Alexander polynomials \\ and HOMFLY polynomials}

The HOMFLY polynomial invariants~\cite{FY} are described by the
skein relationship: 
\be\label{homfly-skein}
a^{-1}H_{+}(a, z)-aH_{-}(a, z)=zH_{O}(a, z). 
 \ee
Comparing (\ref{gA-}) with the HOMFLY skein relationship (\ref{homfly-skein})
we obtain
\be\label{azpq}
 a=q^{1\over4}p^{1\over4}, \quad z=q^{1\over4}p^{-{1\over4}}-q^{-{1\over4}}p^{1\over4}\,.\ee
Substituting this result into (\ref{homfly-skein}), one obtains the generalized Alexander skein relationship (\ref{gA-}).

\vspace{0,5cm}
\noindent
\section{
Generalized Alexander polynomials \\ and Jones polynomials}

The Jones polynomial invariants ~\cite{Jo} can be defined as
\be\label{jones-skein}
t^{-1}V_{+}(t)-tV_{-}(t)=(t^{1\over2}-t^{-{1\over2}})V_{O}(t).\ee 
 Comparing (\ref{gA-}) with the Jones skein relationship (\ref{jones-skein}),
we find that substitution
\be\label{qpt}
q=t^{3},\quad p=t\ee
reduces the generalized Alexander polynomials to Jones ones.

According to results of Section 3, the Jones skein relationship~(\ref{jones-skein})
can be obtained with the help of the proposed three-step algorithm
from $q-$numbers defined as
\be \label{q-jones}
[n]_{q^{3},q}={{\textstyle {q^{3n}-q^{n}}}\over{\textstyle {q^{3}-q}}}\,. \ee

\vspace{0.3cm}
\noindent
\section*{
Acknowledgement}

This work was partially supported by the Special Programme of Division
of Physics and Astronomy of NAS of Ukraine.

\end{document}